# SPARSE RECOVERY IN CONVEX HULLS VIA ENTROPY PENALIZATION[1]

BY VLADIMIR KOLTCHINSKII

*Georgia Institute of Technology*

Let $(X,Y)$ be a random couple in $S \times T$ with unknown distribution $P$ and $(X_1,Y_1),\ldots,(X_n,Y_n)$ be i.i.d. copies of $(X,Y)$. Denote $P_n$ the empirical distribution of $(X_1,Y_1),\ldots,(X_n,Y_n)$. Let $h_1,\ldots,h_N:S \mapsto [-1,1]$ be a dictionary that consists of $N$ functions. For $\lambda \in \mathbb{R}^N$, denote $f_\lambda := \sum_{j=1}^N \lambda_j h_j$. Let $\ell:T \times \mathbb{R} \mapsto \mathbb{R}$ be a given loss function and suppose it is convex with respect to the second variable. Let $(\ell \bullet f)(x,y) := \ell(y; f(x))$. Finally, let $\Lambda \subset \mathbb{R}^N$ be the simplex of all probability distributions on $\{1,\ldots,N\}$. Consider the following penalized empirical risk minimization problem

$$\hat{\lambda}^\varepsilon := \operatorname*{argmin}_{\lambda \in \Lambda}\left[P_n(\ell \bullet f_\lambda) + \varepsilon \sum_{j=1}^N \lambda_j \log \lambda_j\right]$$

along with its distribution dependent version

$$\lambda^\varepsilon := \operatorname*{argmin}_{\lambda \in \Lambda}\left[P(\ell \bullet f_\lambda) + \varepsilon \sum_{j=1}^N \lambda_j \log \lambda_j\right],$$

where $\varepsilon \geq 0$ is a regularization parameter. It is proved that the "approximate sparsity" of $\lambda^\varepsilon$ implies the "approximate sparsity" of $\hat{\lambda}^\varepsilon$ and the impact of "sparsity" on bounding the excess risk of the empirical solution is explored. Similar results are also discussed in the case of entropy penalized density estimation.

**1. Introduction.** Let $S$ and $T$ be measurable spaces with $\sigma$-algebras $\mathcal{S}$ and $\mathcal{T}$, respectively, and let $(X,Y)$ be a random couple in $S \times T$. The distribution of $(X,Y)$ will be denoted by $P$ and the distribution of $X$ by $\Pi$. The training data $(X_1,Y_1),\ldots,(X_n,Y_n)$ consists of $n$ i.i.d. copies of $(X,Y)$ (the distribution $P$ is not known and it is to be estimated based on the data).

Received December 2007; revised May 2008.
[1]Supported in part by NSF Grant DMS-MSPA-06-24841.

We will denote $P_n$ the empirical distribution of the data and will write in what follows

$$Pg = \mathbb{E}g(X,Y) \quad \text{and} \quad P_n g = n^{-1} \sum_{j=1}^{n} g(X_i, Y_i)$$

for functions $g$ on $S \times T$ (as well as for functions on $S$ since they can be also viewed as functions on $S \times T$).

We will be interested in a class of prediction problems in which $Y$ is to be predicted based on an observation of $X$. Prediction rules will be based on the training data $(X_1, Y_1), \ldots, (X_n, Y_n)$.

Let $\ell : T \times \mathbb{R} \mapsto \mathbb{R}_+$ be a loss function. It will be assumed in what follows that, for all $y \in T$, $\ell(y, \cdot)$ is convex. For a function $f : S \mapsto \mathbb{R}$, let $(\ell \bullet f)(x,y) := \ell(y, f(x))$. Then the quantity $P(\ell \bullet f)$ is *the (true) risk* of the prediction rule $f$ and $P_n(\ell \bullet f)$ is the corresponding empirical risk. *The excess risk* of $f$ is defined as

$$\mathcal{E}(f) := P(\ell \bullet f) - \inf_{g:S \mapsto \mathbb{R}} P(\ell \bullet g) = P(\ell \bullet f) - P(\ell \bullet f_*),$$

where the infimum is taken over all measurable functions and it is assumed for simplicity that it is attained at $f_* \in L_2(\Pi)$ (moreover, it will be assumed in what follows that $f_*$ is uniformly bounded by a constant $M$).

Let

$$\mathcal{H} := \{h_1, \ldots, h_N\}$$

be a given finite class of measurable functions from $S$ into $[-1,1]$ called *a dictionary* (of course, it can be assumed instead that the functions in the dictionary are uniformly bounded by an arbitrary constant; the only change will be in the constants in the results below). The dictionary can be an orthonormal system of functions, a union of several orthonormal systems suitable for approximation of the target function $f_*$, a base class of a boosting type algorithm, a set of pretrained estimators in an aggregation problem, etc. Let $\mathcal{P}(\mathcal{H})$ be the set of all probability measures on $\mathcal{H}$. For $\lambda \in \mathcal{P}(\mathcal{H})$, denote $\lambda_j := \lambda(\{h_j\})$ and

$$f_\lambda(x) := \int_{\mathcal{H}} h(x) \lambda(dh) = \sum_{j=1}^{N} \lambda_j h_j(x).$$

Denote $\Lambda := \{(\lambda_1, \ldots, \lambda_N) : \lambda_j \geq 0, j = 1, \ldots, N, \sum_{j=1}^{N} \lambda_j = 1\}$. We will identify (whenever it is convenient) probability measures $\lambda \in \mathcal{P}(\mathcal{H})$ with vectors $(\lambda_1, \ldots, \lambda_N)$ from the simplex $\Lambda$. We will write (with a little abuse of notation) $\lambda = (\lambda_1, \ldots, \lambda_N)$. Clearly, the function $f_\lambda : S \mapsto [-1, 1]$ is a convex combination (a mixture) of functions from the dictionary and the set

$$\mathrm{conv}(\mathcal{H}) := \{f_\lambda : \lambda \in \mathcal{P}(\mathcal{H})\}$$

is the convex hull of $\mathcal{H}$.

As always, define *the entropy* of $\lambda$ as

$$H(\lambda) = -\sum_{j=1}^{N} \lambda_j \log \lambda_j.$$

The Kullback–Leibler divergence between $\lambda, \nu \in \Lambda$ is defined as

$$K(\lambda|\nu) := \sum_{j=1}^{N} \lambda_j \log\left(\frac{\lambda_j}{\nu_j}\right).$$

Denote

$$K(\lambda, \nu) := K(\lambda|\nu) + K(\nu|\lambda).$$

The following penalized empirical risk minimization problem will be studied:

(1.1)
$$\hat{\lambda}^{\varepsilon} := \operatorname*{argmin}_{\lambda \in \mathcal{P}(\mathcal{H})} [P_n(\ell \bullet f_\lambda) - \varepsilon H(\lambda)]$$
$$= \operatorname*{argmin}_{\lambda \in \Lambda} \left[ P_n(\ell \bullet f_\lambda) + \varepsilon \sum_{j=1}^{N} \lambda_j \log \lambda_j \right],$$

where $\varepsilon \geq 0$ is a regularization parameter. Since, for all $y$, $\ell(y, \cdot)$ is convex, the empirical risk $P_n(\ell \bullet f_\lambda)$ is a convex function of $\lambda$. Since also the set $\mathcal{P}(\mathcal{H})$ is convex (it can be identified with the simplex $\Lambda$) and the function $\lambda \mapsto -H(\lambda)$ is convex, this makes the problem (1.1) a convex optimization problem. It is natural to compare this problem with its distribution dependent version

(1.2)
$$\lambda^{\varepsilon} := \operatorname*{argmin}_{\lambda \in \mathcal{P}(\mathcal{H})} [P(\ell \bullet f_\lambda) - \varepsilon H(\lambda)]$$
$$= \operatorname*{argmin}_{\lambda \in \Lambda} \left[ P(\ell \bullet f_\lambda) + \varepsilon \sum_{j=1}^{N} \lambda_j \log \lambda_j \right].$$

In the recent literature, there has been considerable attention to the problem of sparse recovery in a linear span of a given dictionary using penalized empirical risk minimization with $\ell_1$-penalty (this method is called LASSO in the literature on regression), and the current paper is close to this line of work. It has become clear that sparse recovery is possible not always, but only under some geometric assumptions on the dictionary. These assumptions are often described in terms of the properties of the Gram matrix of the dictionary, which in the case of random design models is the matrix

$$H := (\langle h_i, h_j \rangle_{L_2(\Pi)})_{i,j=1,N},$$

and they take form of various conditions on the entries of this matrix ("coherence coefficients"), or on its submatrices (in spirit of "uniform uncertainty principle" or "restricted isometry" conditions). The essence of these assumptions is to try to keep the dictionary not too far from being orthonormal in $L_2(\Pi)$, which in some sense is an ideal case for sparse recovery [see, e.g., Donoho (2006), Candes and Tao (2007), Rudelson and Vershynin (2005), Mendelson, Pajor and Tomczak-Jaegermann (2007), Bunea, Tsybakov and Wegkamp (2007a), van de Geer (2008), Koltchinskii (2008a, 2008b) and Bickel, Ritov and Tsybakov (2008), among many other papers that study both the random design and the fixed design problems].

The idea to use the entropy for complexity regularization is not new in information theory and statistics (recall, e. g., the principle of maximum entropy). In particular, it has been studied recently in connection with the problem of aggregation of statistical estimators by exponential weighting and also in a large number of papers on PAC-Bayesian approach in learning theory [see, e.g., McAllester (1999), Catoni (2004), Audibert (2004), Zhang (2001, 2006a, 2006b) and references therein]. However, we are not aware of any attempt to relate this penalization technique to sparse recovery problems with an exception of a very recent paper by Dalalyan and Tsybakov (2007), where it is done in the context of aggregation with exponential weighting. Moreover, at least at the first glance, the idea of using this type of penalization to achieve sparse recovery seems counterintuitive since the penalty $-H(\lambda)$ attains its minimum at the uniform distribution $\lambda_j = N^{-1}, j = 1, \ldots, N$, and, from this point of view, it penalizes for "sparsity" rather than for "nonsparsity" [in fact, solutions of (1.1), (1.3) can be only "approximately sparse"].

In this paper we follow the approach of Koltchinskii (2005, 2008a), where the problem was studied in the case of $\ell_p$-penalization with $1 \leq p \leq 1 + \frac{c}{\log N}$. This approach is based on separate study of random error $|\mathcal{E}(f_{\hat{\lambda}^\varepsilon}) - \mathcal{E}(f_{\lambda^\varepsilon})|$ and of approximation error $\mathcal{E}(f_{\lambda^\varepsilon})$. It happens that these are two different problems with not entirely the same geometric parameters responsible for the size of each of the two errors, and the geometry of the problem is more subtle than in the standard approach based on conditions on the Gram matrix $H$. In many problems in Statistics and Learning Theory the distribution of the design variable is completely unknown and it is unrealistic to make any restrictive assumptions on its Gram matrix. Because of this reason, it is desirable to study in a more precise way how the excess risk of the solution of (1.1) depends on geometric parameters of the problem.

One of our goals is to show that if $\lambda^\varepsilon$ is "*approximately sparse*" (i.e., this measure is almost concentrated on a small set of atoms), then a similar property is enjoyed by $\hat{\lambda}^\varepsilon$. These sparsity bounds provide a way to control $\|f_{\hat{\lambda}^\varepsilon} - f_{\lambda^\varepsilon}\|_{L_2(\Pi)}$ and $K(\hat{\lambda}^\varepsilon, \lambda^\varepsilon)$ (see Theorems 1 and 2). For instance, we

show that for any set $J \subset \{1, \ldots, N\}$ with $\text{card}(J) = d$ and such that

$$\sum_{j \notin J} \lambda_j^\varepsilon \leq \sqrt{\frac{\log N}{n}},$$

the following bound holds with a high probability:

$$\|f_{\hat{\lambda}^\varepsilon} - f_{\lambda^\varepsilon}\|_{L_2(\Pi)}^2 + \varepsilon K(\hat{\lambda}^\varepsilon; \lambda^\varepsilon) \leq C \frac{d + \log N}{n}.$$

This allow us also to bound "the random error" $|\mathcal{E}(f_{\hat{\lambda}^\varepsilon}) - \mathcal{E}(f_{\lambda^\varepsilon})|$ in terms of "approximate sparsity" of the problem (Theorem 3).

Some further geometric parameters (such as "the alignment coefficient" introduced in the next section) provide a way to control "the approximation error" $\mathcal{E}(f_{\lambda^\varepsilon})$ (see Theorem 4). Namely, suppose there exists a vector $\lambda \in \Lambda$ with the following properties:

(i) $\lambda$ is "*sparse*" [i.e., its support $J = \text{supp}(\lambda)$ is a set of relatively small cardinality];

(ii) the excess risk $\mathcal{E}(f_\lambda)$ is small;

(iii) $\lambda$ is "aligned" nicely with the dictionary (the precise definitions are given in the next section).

Then $\lambda^\varepsilon$ is approximately sparse and its excess risk $\mathcal{E}(f_{\lambda^\varepsilon})$ is small (more precisely, its size is controlled by sparsity of $\lambda$ and its "alignment" with the dictionary). These results ultimately yield oracle inequalities on the excess risk $\mathcal{E}(f_{\hat{\lambda}^\varepsilon})$ showing that this estimation method provides certain degree of adaptation to unknown "sparsity" of the problem (see Corollary 1).

Density estimation problem can be also studied rather naturally in a similar framework. In this problem, the data consists of $n$ independent identically distributed observations $X_1, \ldots, X_n$ in $S$ with common distribution $P$. Suppose that $P$ has density $f_*$ with respect to a $\sigma$-finite measure $\mu$ in $(S, \mathcal{A})$. We will assume that $f_*$ is uniformly bounded by a constant $M$. Let $h_1, \ldots, h_N$ be a large dictionary of probability densities with respect to $\mu$ uniformly bounded by 1 (as in the case of prediction problem discussed above, one can assume that these densities are uniformly bounded by an arbitrary constant resulting in a proper change of constants in the theorems). The goal is to construct an estimator of $f_*$ in the class of mixtures $\{f_\lambda : \lambda \in \Lambda\}$. The underlying assumption is that there exists a "sparse" mixture that approximates the unknown density reasonably well. One can use an estimator based on minimizing the entropy penalized empirical risk with respect to quadratic loss:

$$(1.3) \qquad \hat{\lambda}^\varepsilon := \underset{\lambda \in \Lambda}{\operatorname{argmin}} \left[ \|f_\lambda\|_{L_2(\mu)}^2 - 2P_n f_\lambda + \varepsilon \sum_{j=1}^N \lambda_j \log \lambda_j \right],$$

which is again a convex minimization problem. The corresponding penalized true risk minimization problem is

$$
\begin{aligned}
\lambda^\varepsilon &:= \operatorname*{argmin}_{\lambda \in \Lambda}[\|f_\lambda - f_*\|^2_{L_2(\mu)} - \varepsilon H(\lambda)] \\
&= \operatorname*{argmin}_{\lambda \in \Lambda}\left[\|f_\lambda\|^2_{L_2(\mu)} - 2Pf_\lambda + \varepsilon \sum_{j=1}^N \lambda_j \log \lambda_j\right].
\end{aligned}
\tag{1.4}
$$

Recently, Bunea, Tsybakov and Wegkamp (2007b) studied a similar density estimation problem with $\ell_1$-penalized empirical risk with respect to quadratic loss (and for the linear aggregation instead of convex aggregation). As in the case of prediction problems (regression, classification), we also obtain the bounds characterizing approximate sparsity of the empirical solution in terms of approximate sparsity of the true solution and oracle inequalities for $\|f_{\hat\lambda^\varepsilon} - f_*\|^2_{L_2(\mu)}$ (which is equivalent to considering the excess risk in this problem; see Theorems 5–7, Corollary 2).

## 2. Main results.

2.1. *Assumptions on the loss.* We assume below the following properties of the loss function $\ell$: for all $y \in T$, $\ell(y, \cdot)$ is twice differentiable, $\ell''_u$ is a uniformly bounded function in $T \times \mathbb{R}$ and

$$
\sup_{y \in T} \ell(y; 0) < +\infty, \qquad \sup_{y \in T} |\ell'_u(y; 0)| < +\infty.
$$

Moreover, denote

$$
\tau(R) := \tfrac{1}{2} \inf_{y \in T} \inf_{|u| \leq R} \ell''_u(y, u).
\tag{2.1}
$$

It will be assumed that

$$
\tau(M \vee 1) > 0
$$

(recall that $M$ is a constant such that $\|f_*\|_\infty \leq M$). Without loss of generality, we also assume that $\tau(R) \leq 1, R > 0$ (otherwise, it can be replaced by a lower bound).

There are many important examples of loss functions satisfying these assumptions, most notably, the quadratic loss $\ell(y, u) := (y - u)^2$ in the case when $T \subset \mathbb{R}$ is a bounded set. In this case, $\tau = 1$. In regression problems with a bounded response variable, one can also consider more general loss functions of the form $\ell(y, u) := \phi(y - u)$, where $\phi$ is an even nonnegative convex twice continuously differentiable function with $\phi''$ uniformly bounded in $\mathbb{R}$, $\phi(0) = 0$ and $\phi''(u) > 0$, $u \in \mathbb{R}$. In binary classification setting (i.e., when $T = \{-1, 1\}$), one can choose the loss $\ell(y, u) = \phi(yu)$ with $\phi$ being a nonnegative decreasing convex twice continuously differentiable function such that $\phi''$ is uniformly bounded in $\mathbb{R}$ and $\phi''(u) > 0$, $u \in \mathbb{R}$. The loss function $\phi(u) = \log_2(1 + e^{-u})$ (often called the logit loss) is a typical example.

Note that the condition that the second derivative $\ell''_u$ is uniformly bounded in $T \times \mathbb{R}$ can be replaced by its uniform boundedness in $T \times [-M \vee 1, M \vee 1]$. The constants in the theorems below will then depend on the sup-norm of the second derivative (and, as a consequence, on $M$); otherwise, the results will be the same. This allows one to cover several other popular choices of the loss function, such as the exponential loss $\ell(y, u) := e^{-yu}$ in binary classification.

We will also assume in what follows that $N \geq (\log n)^\gamma$ for some $\gamma > 0$ (this is needed only to avoid additional terms of the order $\frac{\log \log n}{n}$ in several inequalities).

2.2. *Sparsity bounds.* Our first goal is to provide upper bounds on $\|f_{\hat{\lambda}^\varepsilon} - f_{\lambda^\varepsilon}\|_{L_2(\Pi)}$, on $K(\hat{\lambda}^\varepsilon, \lambda^\varepsilon)$ and on $\sum_{j \notin J} \hat{\lambda}^\varepsilon_j$, for an arbitrary subset $J \subset \{1, \ldots, N\}$, in terms of the cardinality of this set $d = \operatorname{card}(J)$ and the measure $\sum_{j \notin J} \lambda^\varepsilon_j$. The idea is to show that if $\lambda^\varepsilon$ is approximately sparse, that is, there exists a small set $J$ such that $\sum_{j \notin J} \lambda^\varepsilon_j$ is also small, then $\hat{\lambda}^\varepsilon$ is approximately sparse, too, with a high probability and the $L_2$-error of approximation of $f_{\lambda^\varepsilon}$ by $f_{\hat{\lambda}^\varepsilon}$ as well as the Kullback–Leibler error of approximation of $\lambda^\varepsilon$ by $\hat{\lambda}^\varepsilon$ are small.

The first result in this direction is the following theorem.

THEOREM 1. *There exist constants $D > 0$ and $C > 0$ depending only on $\ell$ such that, for all $J \subset \{1, \ldots, N\}$ with $d := d(J) = \operatorname{card}(J)$, for all $A \geq 1$ and for all*

$$(2.2) \qquad \varepsilon \geq D\sqrt{\frac{d + A \log N}{n}},$$

*the following bounds hold with probability at least $1 - N^{-A}$:*

$$\sum_{j \notin J} \hat{\lambda}^\varepsilon_j \leq C\left[\sum_{j \notin J} \lambda^\varepsilon_j + \sqrt{\frac{d + A \log N}{n}}\right],$$

$$\sum_{j \notin J} \lambda^\varepsilon_j \leq C\left[\sum_{j \notin J} \hat{\lambda}^\varepsilon_j + \sqrt{\frac{d + A \log N}{n}}\right]$$

*and*

$$\|f_{\hat{\lambda}^\varepsilon} - f_{\lambda^\varepsilon}\|^2_{L_2(\Pi)} + \varepsilon K(\hat{\lambda}^\varepsilon, \lambda^\varepsilon) \leq C\left[\frac{d + A \log N}{n} \vee \sum_{j \notin J} \lambda^\varepsilon_j \sqrt{\frac{d + A \log N}{n}}\right].$$

Note that these bounds hold without any conditions on the dictionary (except the assumption that the functions $h_j$ are uniformly bounded). However, the result is true only for $\varepsilon \geq D\sqrt{\frac{d+A\log N}{n}}$. Since it is not known for which set $J \sum_{j \notin J} \lambda_j^\varepsilon$ is small, it is also not known for which $d$ the condition (2.2) is to be satisfied. In other words, the regularization parameter $\varepsilon$ in this result depends on unknown degree of sparsity of the problem.

In the next theorem, it will be assumed only that $\varepsilon \geq D\sqrt{\frac{A\log N}{n}}$, but there will be more dependence of the bounds on the geometric properties of the dictionary. On the other hand, the error will be controlled not by $d = \operatorname{card}(J)$, but rather by the dimension of a linear space $L$ that provides a good approximation of the functions $\{h_j : j \in J\}$. This dimension can be smaller than $\operatorname{card}(J)$, which makes the bound more precise. Given a subspace $L$ of $L_2(\Pi)$, define

$$U(L) := \sup_{f \in L, \|f\|_{L_2(\Pi)}=1} \|f\|_\infty + 1.$$

It is easy to check that for any $L_2(\Pi)$-orthonormal basis $\phi_1, \ldots, \phi_d$ of $L$,

$$U(L) \leq \max_{1 \leq j \leq d} \|\phi_j\|_\infty \sqrt{d} + 1,$$

where $d := \dim(L)$. In what follows $P_L$ denotes the orthogonal projector onto $L$ and $L^\perp$ denotes the orthogonal complement of $L$. We will be interested in subspaces $L$ for which $\dim(L)$ and $U(L)$ are not very large and, at the same time, functions $\{h_j : j \in J\}$ in the "relevant" part of the dictionary can be approximated well by the functions from $L$ in the sense that the quantity $\max_{j \in J} \|P_{L^\perp} h_j\|_{L_2(\Pi)}$ is relatively small.

THEOREM 2. *Suppose that*

(2.3) $$\varepsilon \geq D\sqrt{\frac{A\log N}{n}}$$

*with a large enough constant $D > 0$ depending only on $\ell$. For all $J \subset \{1, \ldots, N\}$, for all subspaces $L$ of $L_2(\Pi)$ with $d := \dim(L)$ and for all $A \geq 1$, the following bounds hold with probability at least $1 - N^{-A}$ and with a constant $C > 0$ depending only on $\ell$:*

(2.4) $$\sum_{j \notin J} \hat{\lambda}_j^\varepsilon \leq C\left[\sum_{j \notin J} \lambda_j^\varepsilon + \frac{d + A\log N}{n\varepsilon} + \max_{j \in J}\|P_{L^\perp}h_j\|_{L_2(\Pi)} + \frac{U(L)\log N}{n\varepsilon}\right],$$

$$\sum_{j \notin J} \lambda_j^\varepsilon \leq C \left[ \sum_{j \notin J} \hat{\lambda}_j^\varepsilon + \frac{d + A \log N}{n\varepsilon} + \max_{j \in J} \|P_{L^\perp} h_j\|_{L_2(\Pi)} \right.$$

(2.5)
$$\left. + \frac{U(L) \log N}{n\varepsilon} \right]$$

and

$$\|f_{\hat{\lambda}^\varepsilon} - f_{\lambda^\varepsilon}\|_{L_2(\Pi)}^2 + \varepsilon K(\hat{\lambda}^\varepsilon, \lambda^\varepsilon)$$

(2.6)
$$\leq C \left[ \frac{d + A \log N}{n} \vee \sum_{j \notin J} \lambda_j^\varepsilon \sqrt{\frac{A \log N}{n}} \right.$$

$$\left. \vee \max_{j \in J} \|P_{L^\perp} h_j\|_{L_2(\Pi)} \sqrt{\frac{A \log N}{n}} \vee \frac{U(L) \log N}{n} \right].$$

If, for some $J$,

$$\sum_{j \notin J} \lambda_j^\varepsilon \leq \sqrt{\frac{A \log N}{n}}$$

and, for some $L$ with $U(L) \leq \sqrt{d}$, $h_j \in L$, $j \in J$, the bound (2.6) simplifies and becomes

$$\|f_{\hat{\lambda}^\varepsilon} - f_{\lambda^\varepsilon}\|_{L_2(\Pi)}^2 + \varepsilon K(\hat{\lambda}^\varepsilon, \lambda^\varepsilon) \leq C \frac{d + A \log N}{n}.$$

It means that the fact that the dictionary is not orthogonal and even is not linearly independent might actually help to make the random errors $\|f_{\hat{\lambda}^\varepsilon} - f_{\lambda^\varepsilon}\|_{L_2(\Pi)}^2$ and $K(\hat{\lambda}^\varepsilon, \lambda^\varepsilon)$ small: their size is controlled in this case by the dimension $d$ of the linear span $L$ of the "relevant part" of the dictionary $\{h_j : j \in J\}$, and $d$ can be much smaller than $\text{card}(J)$.

2.3. *Random error bounds.* The following result is a simple corollary of Theorems 1, 2 and the properties of the loss function. Denote by $\mathcal{L}$ the linear span of the dictionary $\{h_1, \ldots, h_N\}$ and let $P_\mathcal{L}$ be the orthogonal projector on $\mathcal{L} \subset L_2(P)$. Define

$$g_\varepsilon := P_\mathcal{L}(\ell' \bullet f_{\lambda^\varepsilon}).$$

THEOREM 3. *Under the conditions of Theorem 1, the following bound holds with probability at least $1 - N^{-A}$, with a constant $C > 0$ depending* *only on $\ell$ and with $d = \mathrm{card}(J)$:*

$$|P(\ell \bullet f_{\hat{\lambda}^\varepsilon}) - P(\ell \bullet f_{\lambda^\varepsilon})|$$

(2.7)
$$\leq C\left[\frac{d + A\log N}{n} \vee \sum_{j \notin J} \lambda_j^\varepsilon \sqrt{\frac{d + A\log N}{n}}\right]$$

$$\vee C^{1/2}\|g_\varepsilon\|_{L_2(\Pi)}\left[\frac{d + A\log N}{n} \vee \sum_{j \notin J} \lambda_j^\varepsilon \sqrt{\frac{d + A\log N}{n}}\right]^{1/2}.$$

*Similarly, under the conditions of Theorem 2, with probability at least $1 - N^{-A}$ and with $d = \dim(L)$*

$$|P(\ell \bullet f_{\hat{\lambda}^\varepsilon}) - P(\ell \bullet f_{\lambda^\varepsilon})|$$

$$\leq C\left[\frac{d + A\log N}{n} \vee \left(\sum_{j \notin J} \lambda_j^\varepsilon \vee \max_{j \in J}\|P_{L^\perp}h_j\|_{L_2(\Pi)}\right)\sqrt{\frac{A\log N}{n}}\right.$$

(2.8)
$$\left.\vee \frac{U(L)\log N}{n}\right]$$

$$\vee C^{1/2}\|g_\varepsilon\|_{L_2(\Pi)}\left[\frac{d + A\log N}{n} \vee \left(\sum_{j \notin J} \lambda_j^\varepsilon \vee \max_{j \in J}\|P_{L^\perp}h_j\|_{L_2(\Pi)}\right)\right.$$

$$\left.\sqrt{\frac{A\log N}{n}} \vee \frac{U(L)\log N}{n}\right]^{1/2}.$$

Recall that $f_*$ denotes a function that minimizes the risk $P(\ell \bullet f)$ and it was assumed that $f_*$ is uniformly bounded by a constant $M$. Clearly, by necessary conditions of minimum, we have

$$P(\ell' \bullet f_*)h_j = 0, \qquad j = 1, \ldots, N,$$

so, $\ell' \bullet f_* \in \mathcal{L}^\perp$. Note that for any function $\bar{f}$ uniformly bounded by $M$ and such that $\ell' \bullet \bar{f} \in \mathcal{L}^\perp$ (in particular, for $f_*$) we have

$$\|g_\varepsilon\|_{L_2(\Pi)} = \|P_\mathcal{L}(\ell' \bullet f_{\lambda^\varepsilon})\|_{L_2(P)} = \|P_\mathcal{L}(\ell' \bullet f_{\lambda^\varepsilon} - \ell' \bullet \bar{f})\|_{L_2(P)}$$

$$\leq \|(\ell' \bullet f_{\lambda^\varepsilon} - \ell' \bullet \bar{f})\|_{L_2(P)} \leq C\|f_{\lambda^\varepsilon} - \bar{f}\|_{L_2(\Pi)},$$

where we used the fact that $\ell'$ is Lipschitz with respect to the second variable.

Under the conditions on the loss function, for all $\lambda \in \Lambda$

(2.9) $\quad \mathcal{E}(f_\lambda) \geq \frac{1}{2}\tau(\|f_*\|_\infty \vee 1)\|f_\lambda - f_*\|_{L_2(\Pi)}^2 =: \tau\|f_\lambda - f_*\|_{L_2(\Pi)}^2,$

which easily follows from a version of Taylor expansion for the risk.

To bound the excess risk $\mathcal{E}(f_{\hat{\lambda}^\varepsilon})$, one has to solve two different problems: bounding the random error

$$|\mathcal{E}(f_{\hat{\lambda}^\varepsilon}) - \mathcal{E}(f_{\lambda^\varepsilon})| = |P(\ell \bullet f_{\hat{\lambda}^\varepsilon}) - P(\ell \bullet f_{\lambda^\varepsilon})|$$

and bounding the approximation error $\mathcal{E}(f_{\lambda^\varepsilon})$. Using the above facts, one can easily get from Theorem 3 the following bounds on the random error: under the conditions of Theorem 1, with probability at least $1 - N^{-A}$ and with $d = \mathrm{card}(J)$

$$|\mathcal{E}(f_{\hat{\lambda}^\varepsilon}) - \mathcal{E}(f_{\lambda^\varepsilon})|$$
$$(2.10) \qquad \leq C\left[\frac{d + A\log N}{n} \vee \sum_{j \notin J} \lambda_j^\varepsilon \sqrt{\frac{d + A\log N}{n}}\right]$$
$$\vee C^{1/2}\sqrt{\frac{\mathcal{E}(f_{\lambda^\varepsilon})}{\tau}}\left[\frac{d + A\log N}{n} \vee \sum_{j \notin J} \lambda_j^\varepsilon \sqrt{\frac{d + A\log N}{n}}\right]^{1/2}$$

and under the conditions of Theorem 2, with probability at least $1 - N^{-A}$ and with $d = \dim(L)$

$$|\mathcal{E}(f_{\hat{\lambda}^\varepsilon}) - \mathcal{E}(f_{\lambda^\varepsilon})|$$
$$\leq C\left[\frac{d + A\log N}{n} \vee \left(\sum_{j \notin J} \lambda_j^\varepsilon \vee \max_{j \in J}\|P_{L^\perp} h_j\|_{L_2(\Pi)}\right)\sqrt{\frac{A\log N}{n}}\right.$$
$$(2.11) \qquad\qquad\qquad\qquad\qquad\qquad\qquad \left.\vee \frac{U(L)\log N}{n}\right]$$
$$\vee C^{1/2}\sqrt{\frac{\mathcal{E}(f_{\lambda^\varepsilon})}{\tau}}\left[\frac{d + A\log N}{n} \vee \left(\sum_{j \notin J} \lambda_j^\varepsilon \vee \max_{j \in J}\|P_{L^\perp} h_j\|_{L_2(\Pi)}\right)\right.$$
$$\left.\sqrt{\frac{A\log N}{n}} \vee \frac{U(L)\log N}{n}\right]^{1/2},$$

which reduces the problem to bounding only the approximation error.

2.4. *Approximation error bounds, alignment and oracle inequalities.* To consider the approximation error we need several definitions. For $\lambda \in \Lambda$, denote

$$T_\Lambda(\lambda) := \{v \in \mathbb{R}^N : \exists t > 0 \ \lambda + vt \in \Lambda\}.$$

The set $T_\Lambda(\lambda)$ is the tangent cone of convex set $\Lambda$ at point $\lambda$. Recall that $H$ denotes the Gram matrix of the dictionary in the space $L_2(\Pi)$. Whenever it is convenient, $H$ will be viewed as a linear transformation of $\mathbb{R}^N$. For a vector $w \in \mathbb{R}^N$, let

$$a_H(\Lambda, \lambda, w) := \sup\{\langle w, u\rangle_{\ell_2} : u \in T_\Lambda(\lambda), \|f_u\|_{L_2(\Pi)} = 1\}.$$

We will call this quantity the *a*lignment coefficient of vector $w$, matrix $H$ and convex set $\Lambda$ at point $\lambda \in \Lambda$. Note that

$$\|f_u\|_{L_2(\Pi)}^2 = \langle Hu, u\rangle_{\ell_2} = \langle H^{1/2}u, H^{1/2}u\rangle_{\ell_2}.$$

Therefore, the alignment coefficient can be bounded as follows:

$$a_H(\Lambda, \lambda, w) \leq \sup\{\langle w, u\rangle_{\ell_2} : u \in \mathbb{R}^N, \ \|f_u\|_{L_2(\Pi)} = 1\}$$
$$= \sup_{\|H^{1/2}u\|_{\ell_2}=1} \langle w, u\rangle_{\ell_2} =: \|w\|_H.$$

If $H$ is nonsingular, we can further write

$$\|w\|_H^2 = \sup_{\|H^{1/2}u\|_{\ell_2}=1} \langle H^{-1/2}w, H^{1/2}u\rangle_{\ell_2} = \|H^{-1/2}w\|_{\ell_2}^2.$$

Even when $H$ is singular, we still have

$$\|w\|_H^2 \leq \|H^{-1/2}w\|_{\ell_2}^2,$$

where for $w \in \mathrm{Im}(H^{1/2}) = H^{1/2}\mathbb{R}^N$, one defines

$$\|H^{-1/2}w\|_{\ell_2} := \inf\{\|v\|_{\ell_2} : H^{1/2}v = w\}$$

[which means factorization of the space with respect to $Ker(H^{1/2})$] and for $w \notin \mathrm{Im}(H^{1/2})$ the norm $\|H^{-1/2}w\|_{\ell_2}$ becomes infinite. It is also easy to see that if $J = \mathrm{supp}(w)$, then

$$\|w\|_H \leq \frac{\|w\|_{\ell_2}}{\sqrt{\kappa(J)(1-\rho^2(J))}} \leq \frac{\|w\|_{\ell_\infty}\sqrt{d(J)}}{\sqrt{\kappa(J)(1-\rho^2(J))}},$$

where $d(J) := \mathrm{card}(J)$, $\kappa(J)$ is the minimal eigenvalue of the matrix

$$H_J = (\langle h_i, h_j\rangle_{L_2(\Pi)})_{i,j \in J}$$

and

$$\rho(J) := \sup\left\{\frac{\langle f_1, f_2\rangle_{L_2(\Pi)}}{\|f_1\|_{L_2(\Pi)}\|f_2\|_{L_2(\Pi)}} : f_1 \in L_J, f_2 \in L_{J^c}\right\},$$

$L_J$ denoting the linear span of $\{h_j : j \in J\}$ [see Koltchinskii (2008a), the proof of Proposition 1, for a similar argument]. Measures of linear dependence similar to $\rho(J)$ are known in multivariate statistical analysis as "canonical correlations."

These upper bounds show that the size of the alignment coefficient is controlled by the "sparsity" of the vector $w$ as well as by some characteristics of the dictionary (or its Gram matrix $H$). In particular, for orthonormal dictionaries and for dictionaries that are close enough to being orthonormal [so that $\kappa(J)$ is bounded away from 0 and $\rho^2(J)$ is bounded away from 1], the alignment coefficient is bounded from above by a quantity of the order $\|w\|_{\ell_\infty}\sqrt{d(J)}$. However, the alignment coefficient can be much smaller than this upper bound and it reflects in a more delicate way rather complicated geometric relationships between the vector $w$, the dictionary and the convex set $\Lambda$. Even the quantity $\|H^{-1/2}w\|_{\ell_2}^2$, which is a rough upper bound on the alignment coefficient that does not take into account the geometry of set $\Lambda$, depends not only on the sparsity of $w$, but also on how this vector is aligned with the eigenspaces of $H$. For instance, if $w$ belongs to the linear span of the eigenspaces that correspond only to the eigenvalues of $H$ that are not too small, $\|H^{-1/2}w\|_{\ell_2}^2$ becomes of the order $\|w\|_{\ell_2}^2$. Note also that the geometry of the problem crucially depends on the unknown distribution $\Pi$ of the design variable [since one has to deal with the Hilbert space $L_2(\Pi)$].

For $\lambda \in \mathbb{R}^N$, let $s_j^N(\lambda) := \log(eN^2\lambda_j)$, $j \in \mathrm{supp}(\lambda)$ and $s_j^N(\lambda) := 0$, $j \notin \mathrm{supp}(\lambda)$. Note that, for $j \in \mathrm{supp}(\lambda)$, $s_j^N(\lambda) = \log\lambda_j + 1 + 2\log N$ and $\log\lambda_j + 1$ is the derivative of the function $\lambda \log \lambda$ involved in the definition of the penalty. Let

$$s^N(\lambda) := (s_1^N(\lambda), \ldots, s_N^N(\lambda)).$$

It happens that both the approximation error $\mathcal{E}(f_{\lambda^\varepsilon})$ and the "approximate sparsity" of $\lambda^\varepsilon$ can be controlled by the alignment coefficient of the vector $s_N(\lambda)$ for an arbitrary $\lambda \in \Lambda$. Denote

$$\alpha_N(\lambda) := a_H(\Lambda, \lambda, s^N(\lambda)).$$

THEOREM 4. *There exists a constant $C > 0$ that depends only on $\ell$ and on the constant $M$ (for which $\|f_*\|_\infty \leq M$) such that, for all $\varepsilon > 0$ and all $\lambda \in \Lambda$*

$$(2.12) \qquad \mathcal{E}(f_{\lambda^\varepsilon}) + 2\varepsilon \sum_{j \notin \mathrm{supp}(\lambda)} \lambda_j^\varepsilon \leq 3\mathcal{E}(f_\lambda) + C\bigg(\varepsilon^2 \alpha_N^2(\lambda) + \frac{\varepsilon}{N}\bigg).$$

Theorem 4 and either of the bounds on the random error (2.10) and (2.11) immediately imply oracle inequalities for the excess risk $\mathcal{E}(f_{\hat\lambda^\varepsilon})$. For instance, the next corollary follows from (2.11).

COROLLARY 1. *Under the conditions of Theorem 2, for all $\lambda \in \Lambda$ with $J = \mathrm{supp}(\lambda)$ and for all subspaces $L$ of $L_2(\Pi)$ with $d := \dim(L)$, the following bound holds with probability at least $1 - N^{-A}$ and with a constant $C$* *depending on $\ell$ and on $M$:*

$$\mathcal{E}(f_{\hat{\lambda}^\varepsilon}) \leq 4\mathcal{E}(f_\lambda) + C\bigg(\frac{d + A\log N}{n} + \max_{j\in J}\|P_{L^\perp}h_j\|_{L_2(\Pi)}\sqrt{\frac{A\log N}{n}} + \frac{U(L)\log N}{n} + \varepsilon^2 \alpha_N^2(\lambda) + \frac{\varepsilon}{N}\bigg).$$

2.5. *Density estimation and sparse mixtures recovery.* In the case of density estimation based on entropy penalized empirical risk minimization with quadratic loss, as in (1.3), the results are rather similar to what was described above for prediction problems (regression and classification) and their proofs are quite similar, too.

Recall the notations at the end of the Introduction 1. Recall also the assumptions that the unknown density $f_*$ of distribution $P$ is uniformly bounded by $M$ and the densities in the dictionary $h_j$ are uniformly bounded by 1.

The following results hold.

THEOREM 5. *There exist numerical constants $D > 0$ and $C > 0$ such that, for all $J \subset \{1,\ldots,N\}$ with $d := d(J) = \mathrm{card}(J)$, for all $A \geq 1$ and for all*

$$\varepsilon \geq D\sqrt{\frac{d + A\log N}{n}},$$

*the following bounds hold with probability at least $1 - N^{-A}$:*

$$\sum_{j\notin J}\hat{\lambda}_j^\varepsilon \leq C\bigg[\sum_{j\notin J}\lambda_j^\varepsilon + M^2\sqrt{\frac{d + A\log N}{n}}\bigg],$$

$$\sum_{j\notin J}\lambda_j^\varepsilon \leq C\bigg[\sum_{j\notin J}\hat{\lambda}_j^\varepsilon + M^2\sqrt{\frac{d + A\log N}{n}}\bigg]$$

*and*

$$\|f_{\hat{\lambda}^\varepsilon} - f_{\lambda^\varepsilon}\|_{L_2(\mu)}^2 + \varepsilon K(\hat{\lambda}^\varepsilon,\lambda^\varepsilon) \leq C\bigg[M^2\frac{d + A\log N}{n} \vee \sum_{j\notin J}\lambda_j^\varepsilon\sqrt{\frac{d + A\log N}{n}}\bigg].$$

THEOREM 6. *Suppose that*

$$\varepsilon \geq D\sqrt{\frac{A\log N}{n}}$$

with a large enough numerical constant $D > 0$. For all $J \subset \{1, \ldots, N\}$, for all subspaces $L$ of $L_2(P)$ with $d := \dim(L)$ and for all $A \geq 1$, the following bounds hold with probability at least $1 - N^{-A}$ and with a numerical constant $C > 0$:

(2.13)
$$\sum_{j \notin J} \hat{\lambda}_j^\varepsilon \leq C \left[ \sum_{j \notin J} \lambda_j^\varepsilon + M^2 \frac{d + A \log N}{n\varepsilon} + \max_{j \in J} \|P_{L^\perp} h_j\|_{L_2(P)} \right.$$
$$\left. + \frac{U(L) \log N}{n\varepsilon} \right],$$

(2.14)
$$\sum_{j \notin J} \lambda_j^\varepsilon \leq C \left[ \sum_{j \notin J} \hat{\lambda}_j^\varepsilon + M^2 \frac{d + A \log N}{n\varepsilon} + \max_{j \in J} \|P_{L^\perp} h_j\|_{L_2(P)} \right.$$
$$\left. + \frac{U(L) \log N}{n\varepsilon} \right]$$

and

(2.15)
$$\|f_{\hat{\lambda}^\varepsilon} - f_{\lambda^\varepsilon}\|_{L_2(\mu)}^2 + \varepsilon K(\hat{\lambda}^\varepsilon, \lambda^\varepsilon)$$
$$\leq C \left[ M^2 \frac{d + A \log N}{n} \vee \sum_{j \notin J} \lambda_j^\varepsilon \sqrt{\frac{A \log N}{n}} \right.$$
$$\vee \max_{j \in J} \|P_{L^\perp} h_j\|_{L_2(P)} \sqrt{\frac{A \log N}{n}}$$
$$\left. \vee \frac{U(L) \log N}{n} \right].$$

In the case of density estimation, it makes sense to redefine the alignment coefficient in terms of measure $\mu$:

$$a_H(\Lambda, \lambda, w) := \sup\{\langle w, u \rangle_{\ell_2} : u \in T_\Lambda(\lambda), \|f_u\|_{L_2(\mu)} = 1\},$$
$$\alpha_N(\lambda) := a_H(\Lambda, \lambda, s^N(\lambda)).$$

The approximation error bound then becomes as follows.

THEOREM 7. *There exists a numerical constant $C > 0$ such that, for all $\varepsilon > 0$ and all $\lambda \in \Lambda$*

(2.16)
$$\|f_{\lambda^\varepsilon} - f_*\|_{L_2(\mu)}^2 + 2\varepsilon \sum_{j \notin \mathrm{supp}(\lambda)} \lambda_j^\varepsilon \leq 3\|f_\lambda - f_*\|_{L_2(\mu)}^2$$
$$+ C\left(\varepsilon^2 \alpha_N^2(\lambda) + \frac{\varepsilon}{N}\right).$$

Finally, this results in the following oracle inequality.

COROLLARY 2. *Under the conditions of Theorem 6, for all $\lambda \in \Lambda$ with $J = \mathrm{supp}(\lambda)$ and for all subspaces $L$ of $L_2(\Pi)$ with $d := \dim(L)$, the following bound holds with probability at least $1 - N^{-A}$ and with a numerical constant $C$:*

$$\|f_{\hat{\lambda}^\varepsilon} - f_*\|_{L_2(\mu)}^2 \leq 4\|f_\lambda - f_*\|_{L_2(\mu)}^2$$
$$+ C\bigg(M^2 \frac{d + A\log N}{n} + \max_{j \in J} \|P_{L^\perp} h_j\|_{L_2(\Pi)} \sqrt{\frac{A\log N}{n}}$$
$$+ \frac{U(L)\log N}{n} + \varepsilon^2 \alpha_N^2(\lambda) + \frac{\varepsilon}{N}\bigg).$$

**3. Proofs.** The proofs of Theorems 1 and 2 are quite similar. We give only the proof of Theorem 2.

PROOF OF THEOREM 2. The following necessary conditions of minima in minimization problems defining $\lambda^\varepsilon$ and $\hat{\lambda}^\varepsilon$ will be used to derive sparsity bounds:

$$(3.1) \qquad P(\ell' \bullet f_{\lambda^\varepsilon})(f_{\hat{\lambda}^\varepsilon} - f_{\lambda^\varepsilon}) + \varepsilon \sum_{j=1}^N (\log \lambda_j^\varepsilon + 1)(\hat{\lambda}_j^\varepsilon - \lambda_j^\varepsilon) \geq 0$$

and

$$(3.2) \qquad P_n(\ell' \bullet f_{\hat{\lambda}^\varepsilon})(f_{\hat{\lambda}^\varepsilon} - f_{\lambda^\varepsilon}) + \varepsilon \sum_{j=1}^N (\log \hat{\lambda}_j^\varepsilon + 1)(\hat{\lambda}_j^\varepsilon - \lambda_j^\varepsilon) \leq 0.$$

The inequality (3.1) holds because the directional derivative of the penalized risk function (which is convex)

$$\Lambda \ni \lambda \mapsto P(\ell \bullet f_\lambda) + \varepsilon \sum_{j=1}^N \lambda_j \log \lambda_j$$

at the point of its minimum $\lambda^\varepsilon$ is nonnegative in the direction of any other point of the convex set $\Lambda$, in this case in the direction of $\hat{\lambda}^\varepsilon$. The inequality (3.2) is based on similar considerations in the case of penalized empirical risk (note that in this case the minimum of the convex function is at $\hat{\lambda}^\varepsilon$ and we are differentiating in the direction to the minimal point, not from the minimal point). Subtracting (3.1) from (3.2) and replacing $P$ by $P_n$ in (3.2), we get

$$P((\ell' \bullet f_{\hat{\lambda}^\varepsilon}) - (\ell' \bullet f_{\lambda^\varepsilon}))(f_{\hat{\lambda}^\varepsilon} - f_{\lambda^\varepsilon}) + \varepsilon \sum_{j=1}^{N}(\log \hat{\lambda}_j^\varepsilon - \log \lambda_j^\varepsilon)(\hat{\lambda}_j^\varepsilon - \lambda_j^\varepsilon)$$

(3.3)
$$\leq (P - P_n)(\ell' \bullet f_{\hat{\lambda}^\varepsilon})(f_{\hat{\lambda}^\varepsilon} - f_{\lambda^\varepsilon}).$$

Note that

$$\sum_{j=1}^{N}(\log \hat{\lambda}_j^\varepsilon - \log \lambda_j^\varepsilon)(\hat{\lambda}_j^\varepsilon - \lambda_j^\varepsilon) = \sum_{j=1}^{N}\left(\log \frac{\hat{\lambda}_j^\varepsilon}{\lambda_j^\varepsilon}\right)(\hat{\lambda}_j^\varepsilon - \lambda_j^\varepsilon) = K(\hat{\lambda}^\varepsilon, \lambda^\varepsilon),$$

so bound (3.3) can be also written as

(3.4)
$$P((\ell' \bullet f_{\hat{\lambda}^\varepsilon}) - (\ell' \bullet f_{\lambda^\varepsilon}))(f_{\hat{\lambda}^\varepsilon} - f_{\lambda^\varepsilon}) + \varepsilon K(\hat{\lambda}^\varepsilon; \lambda^\varepsilon)$$
$$\leq (P - P_n)(\ell' \bullet f_{\hat{\lambda}^\varepsilon})(f_{\hat{\lambda}^\varepsilon} - f_{\lambda^\varepsilon}).$$

To extract from this bound some information about approximate sparsity of $\hat{\lambda}^\varepsilon$ note that

(3.5)
$$K(\hat{\lambda}^\varepsilon, \lambda^\varepsilon) = \sum_{j=1}^{N}\left(\log \frac{\hat{\lambda}_j^\varepsilon}{\lambda_j^\varepsilon}\right)(\hat{\lambda}_j^\varepsilon - \lambda_j^\varepsilon)$$
$$\geq \frac{\log 2}{2} \sum_{j:\hat{\lambda}_j^\varepsilon \geq 2\lambda_j^\varepsilon} \hat{\lambda}_j^\varepsilon + \frac{\log 2}{2} \sum_{j:\lambda_j^\varepsilon \geq 2\hat{\lambda}_j^\varepsilon} \lambda_j^\varepsilon.$$

This implies that for any $J \subset \{1, \ldots, N\}$

(3.6)
$$\sum_{j \notin J} \hat{\lambda}_j^\varepsilon \leq 2 \sum_{j \notin J} \lambda_j^\varepsilon + \frac{2}{\log 2} K(\hat{\lambda}^\varepsilon, \lambda^\varepsilon).$$

Similarly,

(3.7)
$$\sum_{j \notin J} \lambda_j^\varepsilon \leq 2 \sum_{j \notin J} \hat{\lambda}_j^\varepsilon + \frac{2}{\log 2} K(\hat{\lambda}^\varepsilon, \lambda^\varepsilon).$$

Therefore, taking into account (3.4),

(3.8)
$$\varepsilon \sum_{j \notin J} \hat{\lambda}_j^\varepsilon \leq 2\varepsilon \sum_{j \notin J} \lambda_j^\varepsilon + \frac{2}{\log 2}(P - P_n)(\ell' \bullet f_{\hat{\lambda}^\varepsilon})(f_{\hat{\lambda}^\varepsilon} - f_{\lambda^\varepsilon}).$$

Since the second derivative of the loss function is bounded away from 0, we also have

$$P((\ell' \bullet f_{\hat{\lambda}^\varepsilon}) - (\ell' \bullet f_{\lambda^\varepsilon}))(f_{\hat{\lambda}^\varepsilon} - f_{\lambda^\varepsilon}) \geq c\|f_{\hat{\lambda}^\varepsilon} - f_{\lambda^\varepsilon}\|^2,$$

where $c = \tau(1)$ (note that $\|f_{\lambda^\varepsilon}\|_\infty \leq 1$ and $\|f_{\hat{\lambda}^\varepsilon}\|_\infty \leq 1$). In view of (3.4), this implies

$$(3.9) \quad c\|f_{\hat{\lambda}^\varepsilon} - f_{\lambda^\varepsilon}\|^2 + \varepsilon K(\hat{\lambda}^\varepsilon, \lambda^\varepsilon) \leq (P - P_n)(\ell' \bullet f_{\hat{\lambda}^\varepsilon})(f_{\hat{\lambda}^\varepsilon} - f_{\lambda^\varepsilon}).$$

Denote

$$\Lambda(\delta; \Delta) := \left\{\lambda \in \Lambda : \|f_\lambda - f_{\lambda^\varepsilon}\|_{L_2(\Pi)} \leq \delta, \sum_{j \notin J} \lambda_j \leq \Delta\right\},$$

$$\alpha_n(\delta; \Delta) := \sup\{|(P_n - P)((\ell' \bullet f_\lambda)(f_\lambda - f_{\lambda^\varepsilon}))| : \lambda \in \Lambda(\delta; \Delta)\}.$$

The following two lemmas are somewhat akin to Lemma 5 in Koltchinskii (2008a). We will give below the proof of Lemma 2 that is needed to complete our proof of Theorem 2. Lemma 1 can be used in a similar way in the proof of Theorem 1, which we skip.

LEMMA 1. *Under the assumptions of Theorem 1, there exists constant $C$ that depends only on $\ell$ such that with probability at least $1 - N^{-A}$, for all*

$$n^{-1/2} \leq \delta \leq 1 \quad \text{and} \quad n^{-1/2} \leq \Delta \leq 1$$

*the following bounds hold:*

$$\alpha_n(\delta; \Delta) \leq \beta_n(\delta; \Delta) := C\left[\delta\sqrt{\frac{d + A\log N}{n}} \vee \Delta\sqrt{\frac{d + A\log N}{n}}\right.$$
$$(3.10)$$
$$\left. \vee \sum_{j \notin J} \lambda_j^\varepsilon \sqrt{\frac{d + A\log N}{n}} \vee \frac{A\log N}{n}\right].$$

LEMMA 2. *Under the assumptions of Theorem 2, there exists constant $C$ that depends only on $\ell$ such that with probability at least $1 - N^{-A}$, for all*

$$(3.11) \quad n^{-1/2} \leq \delta \leq 1 \quad \text{and} \quad n^{-1/2} \leq \Delta \leq 1$$

*the following bounds hold:*

$$\alpha_n(\delta; \Delta) \leq \beta_n(\delta; \Delta)$$
$$(3.12) \quad := C\left[\delta\sqrt{\frac{d + A\log N}{n}} \vee \Delta\sqrt{\frac{A\log N}{n}} \vee \sum_{j \notin J} \lambda_j^\varepsilon \sqrt{\frac{A\log N}{n}}\right.$$
$$\left. \vee \max_{j \in J} \|P_{L^\perp} h_j\|_{L_2(\Pi)} \sqrt{\frac{A\log N}{n}} \vee \frac{U(L)\log N}{n} \vee \frac{A\log N}{n}\right].$$

It follows from Lemma 2 and from (3.8), (3.9) that, for

(3.13) $$\delta = \|f_{\hat{\lambda}^\varepsilon} - f_{\lambda^\varepsilon}\|_{L_2(\Pi)} \qquad \text{and} \qquad \Delta = \sum_{j \notin J} \hat{\lambda}_j^\varepsilon,$$

the following bounds hold with $\beta_n(\delta, \Delta)$ defined in (3.12):

(3.14) $$c\delta^2 \leq \beta_n(\delta, \Delta)$$

and

(3.15) $$\varepsilon\Delta \leq 2\varepsilon \sum_{j \notin J} \lambda_j^\varepsilon + \frac{2}{\log 2}\beta_n(\delta, \Delta)$$

provided that $\delta \geq n^{-1/2}, \Delta \geq n^{-1/2}$. In the case if $\delta < n^{-1/2}$ or $\Delta < n^{-1/2}$ one can replace $\delta$ or $\Delta$, respectively, by $n^{-1/2}$ in the expression for $\beta_n(\delta, \Delta)$ and still have bounds (3.14) and (3.15). The proof below goes through in this case, even with some simplifications. In the main case, when $\delta \geq n^{-1/2}, \Delta \geq n^{-1/2}$, it remains to solve the inequalities (3.14), (3.15) to complete the proof. To this end, note that (3.15) can be rewritten (with a proper adjustment of constant $C$) as

$$\varepsilon\Delta \leq C\Delta\sqrt{\frac{A\log N}{n}}$$
$$+ C\left[\varepsilon \sum_{j \notin J} \lambda_j^\varepsilon \vee \delta\sqrt{\frac{d + A\log N}{n}} \vee \sum_{j \notin J} \lambda_j^\varepsilon\sqrt{\frac{A\log N}{n}}\right.$$
$$\left.\vee \max_{j \in J}\|P_{L^\perp}h_j\|_{L_2(\Pi)}\sqrt{\frac{A\log N}{n}} \vee \frac{U(L)\log N}{n} \vee \frac{A\log N}{n}\right].$$

Under the assumption that the constant $D$ in the condition (2.3) on $\varepsilon$ is larger than 1, the term $\sum_{j \notin J} \lambda_j^\varepsilon \sqrt{\frac{A\log N}{n}}$ in the maximum can be dropped since it smaller than the first term $\varepsilon \sum_{j \notin J} \lambda_j^\varepsilon$. If $D \geq 2C$, the bound can be further rewritten as

$$\varepsilon\Delta \leq C\left[\varepsilon \sum_{j \notin J} \lambda_j^\varepsilon \vee \delta\sqrt{\frac{d + A\log N}{n}}\right.$$
$$\vee \max_{j \in J}\|P_{L^\perp}h_j\|_{L_2(\Pi)}\sqrt{\frac{A\log N}{n}}$$
$$\left.\vee \frac{U(L)\log N}{n} \vee \frac{A\log N}{n}\right]$$

(again with an adjustment of $C$). To get a bound on $\Delta$, it is enough to solve the inequality separately for each term in the maximum and take the maximum of the solutions. This yields

$$\Delta \leq C \Bigg[\sum_{j \notin J} \lambda_j^\varepsilon \vee \frac{\delta}{\varepsilon}\sqrt{\frac{d + A \log N}{n}}$$

$$\vee \max_{j \in J} \|P_{L^\perp} h_j\|_{L_2(\Pi)} \frac{1}{\varepsilon}\sqrt{\frac{A \log N}{n}} \vee \frac{U(L)\log N}{n\varepsilon} \vee \frac{A \log N}{n\varepsilon}\Bigg].$$

Under the assumption (2.3) on $\varepsilon$ (with $D \geq 1$), this can be further simplified and the bound becomes

$$\Delta \leq \Delta(\delta) := C \Bigg[\sum_{j \notin J} \lambda_j^\varepsilon \vee \frac{\delta}{\varepsilon}\sqrt{\frac{d + A \log N}{n}}$$

$$\vee \max_{j \in J} \|P_{L^\perp} h_j\|_{L_2(\Pi)} \vee \frac{U(L)\log N}{n\varepsilon} \vee \sqrt{\frac{A \log N}{n}}\Bigg].$$

Let us now substitute $\Delta(\delta)$ instead of $\Delta$ in (3.14) [note than $\beta_n(\delta,\Delta)$ is nondecreasing in $\Delta$]. This easily gives the following bound on $\delta$:

$$\delta^2 \leq C\Bigg[\delta\sqrt{\frac{d + A \log N}{n}} \vee \frac{\delta}{\varepsilon}\sqrt{\frac{d + A \log N}{n}}\sqrt{\frac{A \log N}{n}}$$

$$\vee \frac{U(L)\log N}{n\varepsilon}\sqrt{\frac{A \log N}{n}} \vee \sum_{j \notin J} \lambda_j^\varepsilon \sqrt{\frac{A \log N}{n}}$$

$$\vee \max_{j \in J} \|P_{L^\perp} h_j\|_{L_2(\Pi)} \sqrt{\frac{A \log N}{n}} \vee \frac{U(L)\log N}{n} \vee \frac{A \log N}{n}\Bigg],$$

and the second and the third terms in the maximum can be dropped again since $\frac{1}{\varepsilon}\sqrt{\frac{A \log N}{n}} \leq 1$. Thus, we have

$$\delta^2 \leq C\Bigg[\delta\sqrt{\frac{d + A \log N}{n}} \vee \sum_{j \notin J} \lambda_j^\varepsilon \sqrt{\frac{A \log N}{n}}$$

$$\vee \max_{j \in J} \|P_{L^\perp} h_j\|_{L_2(\Pi)}\sqrt{\frac{A \log N}{n}} \vee \frac{U(L)\log N}{n} \vee \frac{A \log N}{n}\Bigg],$$

which gives the following bound on $\delta^2$:

$$\delta^2 \leq C\Bigg[\frac{d + A \log N}{n} \vee \sum_{j \notin J} \lambda_j^\varepsilon \sqrt{\frac{A \log N}{n}}$$

SPARSE RECOVERY    21(3.16)
$$\vee \max_{j\in J}\|P_{L^\perp}h_j\|_{L_2(\Pi)}\sqrt{\frac{A\log N}{n}} \vee \frac{U(L)\log N}{n}\Bigg].$$

This can be substituted back into the expression for $\Delta(\delta)$ yielding the bound on $\Delta$:

$$\Delta \leq C\Bigg[\sum_{j\notin J}\lambda_j^\varepsilon \vee \frac{d+A\log N}{n\varepsilon} \vee \Big(\sum_{j\notin J}\lambda_j^\varepsilon\Big)^{1/2}\frac{1}{\varepsilon}\Big(\frac{A\log N}{n}\Big)^{1/4}\sqrt{\frac{d+A\log N}{n}}$$

$$\vee \sqrt{\frac{U(L)\log N}{n\varepsilon}}\sqrt{\frac{d+A\log N}{n}}$$

$$\vee \max_{j\in J}\|P_{L^\perp}h_j\|_{L_2(\Pi)}^{1/2}\frac{1}{\varepsilon}\Big(\frac{A\log N}{n}\Big)^{1/4}\sqrt{\frac{d+A\log N}{n}}$$

$$\vee \max_{j\in J}\|P_{L^\perp}h_j\|_{L_2(\Pi)} \vee \frac{U(L)\log N}{n\varepsilon} \vee \sqrt{\frac{A\log N}{n}}\Bigg],$$

which, using the inequality $ab\leq (a^2+b^2)/2$ and the condition $\frac{1}{\varepsilon}\sqrt{\frac{A\log N}{n}}\leq 1$, can be simplified and rewritten as

(3.17)
$$\Delta \leq C\Bigg[\sum_{j\notin J}\lambda_j^\varepsilon \vee \frac{d+A\log N}{n\varepsilon} \vee \max_{j\in J}\|P_{L^\perp}h_j\|_{L_2(\Pi)}$$
$$\vee \frac{U(L)\log N}{n\varepsilon} \vee \sqrt{\frac{A\log N}{n}}\Bigg]$$

with a proper change of $C$ (still depending only on $\ell$). Now we can substitute (3.16) and (3.17) in the expression for $\beta_n(\delta,\Delta)$. We skip the details that are simple and similar to the bounds earlier in the proof. In view of Lemma 2, this gives the following bound on $\alpha_n(\delta,\Delta)$ that holds for $\delta,\Delta$ defined by (3.13) with probability at least $1-N^{-A}$:

$$\alpha_n(\delta,\Delta)\leq C\Bigg[\frac{d+A\log N}{n}+\sum_{j\notin J}\lambda_j^\varepsilon\sqrt{\frac{A\log N}{n}}$$

$$\vee \max_{j\in J}\|P_{L^\perp}h_j\|_{L_2(\Pi)}\sqrt{\frac{A\log N}{n}} \vee \frac{U(L)\log N}{n}\Bigg].$$

Together with (3.9) this yields the bound

$$c\|f_{\hat\lambda^\varepsilon}-f_{\lambda^\varepsilon}\|_{L_2(\Pi)}^2+\varepsilon K(\hat\lambda^\varepsilon,\lambda^\varepsilon)$$

$$\text{(3.18)} \quad \leq C\left[\frac{d + A\log N}{n} + \sum_{j \notin J} \lambda_j^\varepsilon \sqrt{\frac{A\log N}{n}}\right.$$
$$\left. \vee \max_{j \in J} \|P_{L^\perp} h_j\|_{L_2(\Pi)} \sqrt{\frac{A\log N}{n}} \vee \frac{U(L)\log N}{n}\right],$$

which is equivalent to (2.6). Bound (2.4) follows immediately from bound (3.17) (under the assumption on $\varepsilon$, the term $\sqrt{\frac{A\log N}{n}}$ is smaller than $\frac{d + A\log N}{n\varepsilon}$, so, it can be discarded), and bound (2.5) follows from (3.7) and (3.18), which completes the proof. □

PROOF OF LEMMA 2. The proof relies on Talagrand's concentration inequality for empirical processes as well as on Rademacher symmetrization and contraction inequalities [see, e.g., Koltchinskii (2006) or Massart (2007) for their formulations in a form convenient for our purposes]. By Talagrand's concentration inequality, with probability at least $1 - e^{-t}$

$$\text{(3.19)} \quad \alpha_n(\delta; \Delta) \leq 2\left[\mathbb{E}\alpha_n(\delta; \Delta) + C\delta\sqrt{\frac{t}{n}} + \frac{Ct}{n}\right]$$

and, by symmetrization inequality,

$$\mathbb{E}\alpha_n(\delta; \Delta) \leq 2\mathbb{E}\sup\{|R_n((\ell' \bullet f_\lambda)(f_\lambda - f_{\lambda^\varepsilon}))| : \lambda \in \Lambda(\delta; \Delta)\}.$$

Since

$$\ell'(f_\lambda(\cdot))(f_\lambda(\cdot) - f_{\lambda^\varepsilon}(\cdot)) = \ell'(f_{\lambda^\varepsilon}(\cdot) + u)u|_{u = f_\lambda(\cdot) - f_{\lambda^\varepsilon}(\cdot)}$$

and the function

$$[-1, 1] \ni u \mapsto \ell'(f_{\lambda^\varepsilon}(\cdot) + u)u$$

is Lipschitz with a constant $C$ depending only on $\ell$, the application of Rademacher contraction inequality yields the bound

$$\text{(3.20)} \quad \mathbb{E}\alpha_n(\delta; \Delta) \leq C\mathbb{E}\sup\{|R_n(f_\lambda - f_{\lambda^\varepsilon})| : \lambda \in \Lambda(\delta; \Delta)\}.$$

Now we use the following representation

$$\text{(3.21)} \quad \begin{aligned} f_\lambda - f_{\lambda^\varepsilon} &= P_L(f_\lambda - f_{\lambda^\varepsilon}) + \sum_{j \in J}(\lambda_j - \lambda_j^\varepsilon)P_{L^\perp} h_j \\ &\quad + \sum_{j \notin J}(\lambda_j - \lambda_j^\varepsilon)P_{L^\perp} h_j. \end{aligned}$$

Clearly, for all $\lambda \in \Lambda(\delta, \Delta)$,

$$\|P_L(f_\lambda - f_{\lambda^\varepsilon})\|_{L_2(\Pi)} \leq \|f_\lambda - f_{\lambda^\varepsilon}\|_{L_2(\Pi)} \leq \delta$$

and $P_L(f_\lambda - f_{\lambda^\varepsilon}) \in L$, which is a $d$-dimensional subspace. Therefore,

$$\mathbb{E}\sup\{|R_n(P_L(f_\lambda - f_{\lambda^\varepsilon}))| : \lambda \in \Lambda(\delta;\Delta)\} \le C\delta\sqrt{\frac{d}{n}}$$

[see, e.g., Koltchinskii (2006), Section 2, Example 1]. On the other hand, since $\lambda, \lambda^\varepsilon \in \Lambda$, we have $\sum_{j\in J} |\lambda_j - \lambda_j^\varepsilon| \le 2$ and

$$\mathbb{E}\sup\left\{\left|R_n\left(\sum_{j\in J}(\lambda_j - \lambda_j^\varepsilon)P_{L^\perp}h_j\right)\right| : \lambda \in \Lambda(\delta;\Delta)\right\} \le 2\mathbb{E}\max_{j\in J}|R_n(P_{L^\perp}h_j)|.$$

We now proceed with rather well-known approach to bounding the sup-norm of Rademacher sums:

$$\mathbb{E}\max_{j\in J}|R_n(P_{L^\perp}h_j)| \le C\mathbb{E}\max_{j\in J}\|P_{L^\perp}h_j\|_{L_2(\Pi_n)}\sqrt{\frac{\log\operatorname{card}(J)}{n}}$$

$$\le C\max_{j\in J}\|P_{L^\perp}h_j\|_{L_2(\Pi)}\sqrt{\frac{\log N}{n}}$$

$$+ \sqrt{\mathbb{E}\max_{j\in J}|\|P_{L^\perp}h_j\|_{L_2(\Pi_n)}^2 - \|P_{L^\perp}h_j\|_{L_2(\Pi)}^2|}\sqrt{\frac{\log N}{n}}.$$

Note that

$$\|P_{L^\perp}h_j\|_\infty \le \|P_L h_j\|_\infty + \|h_j\|_\infty \le (U(L)-1)\|P_L h_j\|_{L_2(\Pi)} + 1$$
$$\le (U(L)-1)\|h_j\|_{L_2(\Pi)} + 1 \le U(L).$$

We use symmetrization inequality together with Rademacher contraction inequality to get the following bound

$$\mathbb{E}\max_{j\in J}|R_n(P_{L^\perp}h_j)|$$

$$\le C\left[\max_{j\in J}\|P_{L^\perp}h_j\|_{L_2(\Pi)}\sqrt{\frac{\log N}{n}}\right.$$

$$\left.+ \sqrt{\max_{j\in J}\|P_{L^\perp}h_j\|_\infty \mathbb{E}\max_{j\in J}|R_n(P_{L^\perp}h_j)|}\sqrt{\frac{\log N}{n}}\right].$$

The last inequality can be solved for

$$\mathbb{E}\max_{j\in J}|R_n(P_{L^\perp}h_j)|,$$

which gives the bound

$$\mathbb{E}\max_{j\in J}|R_n(P_{L^\perp}h_j)| \le C\left[\max_{j\in J}\|P_{L^\perp}h_j\|_{L_2(\Pi)}\sqrt{\frac{\log N}{n}} + U(L)\frac{\log N}{n}\right].$$

Quite similarly, we have for all $\lambda \in \Lambda(\delta, \Delta)$

$$\sum_{j \notin J} |\lambda_j - \lambda_j^{\varepsilon}| \leq \Delta + \sum_{j \notin J} \lambda_j^{\varepsilon}$$

and

$$\mathbb{E} \sup \left\{ \left| R_n \left( \sum_{j \notin J} (\lambda_j - \lambda_j^{\varepsilon}) P_{L^{\perp}} h_j \right) \right| : \lambda \in \Lambda(\delta; \Delta) \right\}$$

$$\leq \left( \Delta + \sum_{j \notin J} \lambda_j^{\varepsilon} \right) \mathbb{E} \max_{j \notin J} |R_n(P_{L^{\perp}} h_j)|.$$

Repeating what we have done in the case of $j \in J$, we get

$$\mathbb{E} \max_{j \notin J} |R_n(P_{L^{\perp}} h_j)| \leq C \left[ \sqrt{\frac{\log N}{n}} + U(L) \frac{\log N}{n} \right],$$

where we used the fact that

$$\|P_{L^{\perp}} h_j\|_{L_2(\Pi)} \leq \|h_j\|_{L_2(\Pi)} \leq 1.$$

It remains to recall representation (3.21) and bound (3.20) to show that

$$\mathbb{E}\alpha_n(\delta, \Delta)$$

$$\leq C \left[ \delta \sqrt{\frac{d}{n}} \vee \Delta \sqrt{\frac{\log N}{n}} \vee \sum_{j \notin J} \lambda_j^{\varepsilon} \sqrt{\frac{\log N}{n}} \right.$$

(3.22)

$$\vee \max_{j \in J} \|P_{L^{\perp}} h_j\|_{L_2(\Pi)} \sqrt{\frac{\log N}{n}} \vee \Delta U(L) \frac{\log N}{n}$$

$$\left. \vee \sum_{j \notin J} \lambda_j^{\varepsilon} U(L) \frac{\log N}{n} \vee U(L) \frac{\log N}{n} \right],$$

which can be bounded further as

$$\mathbb{E}\alpha_n(\delta, \Delta) \leq C \left[ \delta \sqrt{\frac{d}{n}} \vee \Delta \sqrt{\frac{\log N}{n}} \vee \sum_{j \notin J} \lambda_j^{\varepsilon} \sqrt{\frac{\log N}{n}} \right.$$

(3.23)

$$\left. \vee \max_{j \in J} \|P_{L^{\perp}} h_j\|_{L_2(\Pi)} \sqrt{\frac{\log N}{n}} \vee \frac{U(L) \log N}{n} \right].$$

This can be plugged in (3.19) to get that with probability $1 - e^{-t}$

$$\alpha_n(\delta, \Delta) \leq \tilde{\beta}_n(\delta, \Delta, t)$$

$$(3.24) \quad := C\left[\delta\sqrt{\frac{d}{n}} \vee \Delta\sqrt{\frac{\log N}{n}} \vee \sum_{j\notin J}\lambda_j^\varepsilon\sqrt{\frac{\log N}{n}}\right.$$

$$\left.\vee \max_{j\in J}\|P_{L^\perp}h_j\|_{L_2(\Pi)}\sqrt{\frac{\log N}{n}} \vee \frac{U(L)\log N}{n} \vee \delta\sqrt{\frac{t}{n}} \vee \frac{t}{n}\right]$$

with a constant $C > 0$ depending only on $\ell$.

We will make the above bound uniform in $\delta, \Delta$ satisfying (3.11). To this end, define

$$\delta_j := 2^{-j} \quad \text{and} \quad \Delta_j := 2^{-j}$$

and replace $t$ by $t + 2\log(j+1) + 2\log(k+1)$. The union bound implies that with probability at least

$$1 - \sum_{j,k\geq 0}\exp\{-t - 2\log(j+1) - 2\log(k+1)\}$$

$$= 1 - \left(\sum_{j\geq 0}(j+1)^{-2}\right)^2 \exp\{-t\} \geq 1 - 4e^{-t},$$

for all $\delta$ and $\Delta$ satisfying (3.11), and for $j, k$ such that

$$\delta \in (\delta_{j+1}, \delta_j] \quad \text{and} \quad \Delta \in (\Delta_{k+1}, \Delta_k],$$

the following bound holds:

$$\alpha_n(\delta; \Delta) \leq \tilde{\beta}_n(\delta_j, \Delta_k, t + 2\log j + 2\log k).$$

Since

$$2\log j \leq 2\log\log_2\left(\frac{1}{\delta_j}\right) \leq 2\log\log_2\left(\frac{2}{\delta}\right)$$

and

$$2\log k \leq 2\log\log_2\left(\frac{2}{\Delta}\right),$$

we have

$$\tilde{\beta}_n(\delta_j, \Delta_k, t + 2\log j + 2\log k)$$
$$\leq \tilde{\beta}_n\left(2\delta, 2\Delta, t + 2\log\log_2\left(\frac{2}{\delta}\right) + 2\log\log_2\left(\frac{2}{\Delta}\right)\right) =: \bar{\beta}_n(\delta; \Delta; t)$$

and, therefore, with probability at least $1 - 4e^{-t}$, for all $\delta$ and $\Delta$ satisfying (3.11),

$$\alpha_n(\delta; \Delta) \leq \bar{\beta}_n(\delta; \Delta; t).$$

Let $t = A \log N + \log 4$ (so that $4e^{-t} = N^{-A}$). Then, with some constant $C$ that depends only on $\ell$,

$$\bar{\beta}_n(\delta; \Delta; t)$$
$$\leq C \bigg[ \delta \sqrt{\frac{d}{n}} \vee \delta \sqrt{\frac{A \log N}{n}} \vee \delta \sqrt{\frac{2 \log \log_2(2/\delta)}{n}} \vee \delta \sqrt{\frac{2 \log \log_2(2/\Delta)}{n}}$$
$$\vee \Delta \sqrt{\frac{\log N}{n}} \vee \sum_{j \notin J} \lambda_j^\varepsilon \sqrt{\frac{\log N}{n}} \vee \max_{j \in J} \|P_{L^\perp} h_j\|_{L_2(\Pi)} \sqrt{\frac{\log N}{n}}$$
$$\vee \frac{U(L) \log N}{n} \vee \frac{2 \log \log_2(2/\delta)}{n} \vee \frac{2 \log \log_2(2/\Delta)}{n} \vee \frac{A \log N}{n} \bigg].$$

For all $\delta$ and $\Delta$ satisfying (3.11),

$$\frac{2 \log \log_2(2/\delta)}{n} \leq C \frac{\log \log n}{n}$$

and

$$\frac{2 \log \log_2(2/\Delta)}{n} \leq C \frac{\log \log n}{n}.$$

By assumptions on $N, n$, $A \log N \geq \gamma \log \log n$. Therefore, for $\delta$ and $\Delta$ satisfying (3.11),

$$\alpha_n(\delta, \Delta) \leq \bar{\beta}_n(\delta; \Delta; t) \leq C \bigg[ \delta \sqrt{\frac{d}{n}} \vee \delta \sqrt{\frac{A \log N}{n}} \vee \Delta \sqrt{\frac{\log N}{n}}$$

(3.25)
$$\vee \sum_{j \notin J} \lambda_j^\varepsilon \sqrt{\frac{\log N}{n}} \vee \max_{j \in J} \|P_{L^\perp} h_j\|_{L_2(\Pi)} \sqrt{\frac{\log N}{n}}$$
$$\vee \frac{U(L) \log N}{n} \vee \frac{A \log N}{n} \bigg],$$

which holds with probability at least $1 - N^{-A}$. □

PROOF OF THEOREM 3. The proof easily follows from the fact that under the conditions on the loss function we have

$$(\ell \bullet f_{\hat{\lambda}^\varepsilon})(x, y) - (\ell \bullet f_{\lambda^\varepsilon})(x, y) = (\ell' \bullet f_{\lambda^\varepsilon})(x, y)(f_{\hat{\lambda}^\varepsilon} - f_{\lambda^\varepsilon})(x) + R(x, y),$$

where

$$|R(x, y)| \leq C(f_{\hat{\lambda}^\varepsilon} - f_{\lambda^\varepsilon})^2(x).$$

Integrating with respect to $P$ yields

$$|P(\ell \bullet f_{\hat{\lambda}^\varepsilon}) - P(\ell \bullet f_{\lambda^\varepsilon}) - P(\ell' \bullet f_{\lambda^\varepsilon})(f_{\hat{\lambda}^\varepsilon} - f_{\lambda^\varepsilon})| \leq C \|f_{\hat{\lambda}^\varepsilon} - f_{\lambda^\varepsilon}\|_{L_2(\Pi)}^2.$$

It is enough to observe that

$$|P(\ell' \bullet f_{\lambda^\varepsilon})(f_{\hat{\lambda}^\varepsilon} - f_{\lambda^\varepsilon})|$$
$$= |\langle \ell' \bullet f_{\lambda^\varepsilon}, f_{\hat{\lambda}^\varepsilon} - f_{\lambda^\varepsilon} \rangle_{L_2(P)}| = |\langle P_{\mathcal{L}}(\ell' \bullet f_{\lambda^\varepsilon}), f_{\hat{\lambda}^\varepsilon} - f_{\lambda^\varepsilon} \rangle_{L_2(P)}|$$
$$\leq \|g_\varepsilon\|_{L_2(P)} \|f_{\hat{\lambda}^\varepsilon} - f_{\lambda^\varepsilon}\|_{L_2(\Pi)}$$

and to use the bounds of Theorems 1 and 2. $\square$

PROOF OF THEOREM 4. The following bound immediately follows from the definition of $\lambda^\varepsilon$: for all $\lambda \in \Lambda$,

$$\mathcal{E}(f_{\lambda^\varepsilon}) + \varepsilon \sum_{j=1}^N \lambda_j^\varepsilon \log(N^2 \lambda_j^\varepsilon) \leq \mathcal{E}(f_\lambda) + \varepsilon \sum_{j=1}^N \lambda_j \log(N^2 \lambda_j).$$

Denoting $J_\lambda = \operatorname{supp}(\lambda)$ and using the convexity of the function $u \mapsto u \log(N^2 u)$ and the fact that its derivative is $\log(eN^2 u)$, we get

$$\mathcal{E}(f_{\lambda^\varepsilon}) + \varepsilon \sum_{j \notin J_\lambda} \lambda_j^\varepsilon \log(N^2 \lambda_j^\varepsilon) \leq \mathcal{E}(f_\lambda) + \varepsilon \sum_{j \in J_\lambda} (\lambda_j \log(N^2 \lambda_j) - \lambda_j^\varepsilon \log(N^2 \lambda_j^\varepsilon))$$
$$\leq \mathcal{E}(f_\lambda) + \varepsilon \sum_{j \in J_\lambda} \log(eN^2 \lambda_j)(\lambda_j - \lambda_j^\varepsilon),$$

which, by the definition of $\alpha_N(\lambda)$, can be further bounded by

$$\mathcal{E}(f_\lambda) + \varepsilon |\alpha_N(\lambda)| \|f_\lambda - f_{\lambda^\varepsilon}\|_{L_2(\Pi)}.$$

Next we use obvious bounds [recall (2.9)]

$$\|f_\lambda - f_{\lambda^\varepsilon}\|_{L_2(\Pi)} \leq \|f_\lambda - f_*\|_{L_2(\Pi)} + \|f_{\lambda^\varepsilon} - f_*\|_{L_2(\Pi)} \leq \sqrt{\frac{\mathcal{E}(f_\lambda)}{\tau}} + \sqrt{\frac{\mathcal{E}(f_{\lambda^\varepsilon})}{\tau}}$$

to get

$$\mathcal{E}(f_{\lambda^\varepsilon}) + \varepsilon \sum_{j \notin J_\lambda} \lambda_j^\varepsilon \log(N^2 \lambda_j^\varepsilon) \leq \mathcal{E}(f_\lambda) + \varepsilon |\alpha_N(\lambda)| \left( \sqrt{\frac{\mathcal{E}(f_\lambda)}{\tau}} + \sqrt{\frac{\mathcal{E}(f_{\lambda^\varepsilon})}{\tau}} \right).$$

Since

$$\varepsilon |\alpha_N(\lambda)| \sqrt{\frac{\mathcal{E}(f_{\lambda^\varepsilon})}{\tau}} \leq \frac{1}{2} \frac{\varepsilon^2 \alpha_N^2(\lambda)}{\tau} + \frac{1}{2} \mathcal{E}(f_{\lambda^\varepsilon})$$

and

$$\varepsilon |\alpha_N(\lambda)| \sqrt{\frac{\mathcal{E}(f_\lambda)}{\tau}} \leq \frac{1}{2} \frac{\varepsilon^2 \alpha_N^2(\lambda)}{\tau} + \frac{1}{2} \mathcal{E}(f_\lambda),$$

this yields

$$\frac{1}{2}\mathcal{E}(f_{\lambda^\varepsilon}) + \varepsilon \sum_{j \notin J_\lambda} \lambda_j^\varepsilon \log(N^2 \lambda_j^\varepsilon) \leq \frac{3}{2}\mathcal{E}(f_\lambda) + \frac{\varepsilon^2 \alpha_N^2(\lambda)}{\tau}.$$

Note that also

$$\sum_{j \notin J_\lambda} \lambda_j^\varepsilon \log(N^2 \lambda_j^\varepsilon) = \sum_{j \notin J_\lambda : \lambda_j^\varepsilon \geq eN^{-2}} \lambda_j^\varepsilon \log(N^2 \lambda_j^\varepsilon) + \sum_{j \notin J_\lambda : \lambda_j^\varepsilon < eN^{-2}} \lambda_j^\varepsilon \log(N^2 \lambda_j^\varepsilon)$$

$$\geq \sum_{j \notin J_\lambda : \lambda_j^\varepsilon \geq eN^{-2}} \lambda_j^\varepsilon - \frac{1}{eN},$$

where we used the fact that the function $t \mapsto t \log(N^2 t)$ is bounded from below by $-\frac{1}{eN^2}$. Thus,

$$\sum_{j \notin J_\lambda} \lambda_j^\varepsilon \log(N^2 \lambda_j^\varepsilon) \geq \sum_{j \notin J_\lambda} \lambda_j^\varepsilon - \sum_{j \notin J_\lambda : \lambda_j^\varepsilon < eN^{-2}} \lambda_j^\varepsilon - \frac{1}{eN} \geq \sum_{j \notin J_\lambda} \lambda_j^\varepsilon - (e + e^{-1})\frac{1}{N}.$$

Therefore, we get

$$\mathcal{E}(f_{\lambda^\varepsilon}) + 2\varepsilon \sum_{j \notin J_\lambda} \lambda_j^\varepsilon \leq 3\mathcal{E}(f_\lambda) + \frac{2\varepsilon^2 \alpha_N^2(\lambda)}{\tau} + 2(e + e^{-1})\frac{\varepsilon}{N},$$

which implies the result. □

The results concerning penalized density estimation can be proved quite similarly.

**Acknowledgment.** The author is thankful to anonymous referee for making several helpful suggestions, in particular, for pointing out a slight improvement of the last term in the bounds of Theorem 4 and Corollary 1.

SCHOOL OF MATHEMATICS
GEORGIA INSTITUTE OF TECHNOLOGY
ATLANTA, GEORGIA 30332-0160
USA
E-MAIL: vlad@math.gatech.edu